# Spectral Analysis and Stability of Wave Equations with Dispersive Nonlinearity

[1]Umar Muhammad Dauda and [2]Lawal Ja'afaru
umarmd2021@gmail.com
[1,2]Department of Mathematics, Faculty of Computing & Mathematical Sciences, Kano University of Science and Technology, Wudil, Kano, Nigeria.

**Abstract:** This study employs spectral methods to capture the behaviour of wave equation with dispersive-nonlinearity. We describe the evolution of hump initial data and track the conservation of the mass and energy functionals. The dispersive-nonlinearity results to solution in an extended Schwartz space via analytic approach. We construct numerical schemes based on spectral methods to simulate soliton interactions under Schwartzian initial data. The computational analysis includes validation of energy and mass conservation to ensure numerical accuracy. Results show that initial data from the Schwartz space decompose into smaller wave-packets due to the weaker dispersive-nonlinearity but leads to wave collapse as a result of stronger dispersive-nonlinearity. We conjecture that the hyperbolic equation with a positive nonlinearity and exponent $\sigma \geq 2$ admits global solutions, while lower exponents lead to localized solutions. A stability analysis of solitonic solutions of the equation is provided via the perturbation approach.

## Introduction

Nonlinearity is an interesting phenomenon, which is, apparently, everywhere. However, dispersion, that is ubiquitous in nature, when it interacts with nonlinearity, it gives more interesting behaviour. These interactions could result in the the so-called solitons. Nonlinear dispersive equations appear in many fields such as water-wave models, nonlinear optics and their likes. Example of such equations include the Korteweg-de-Vries (KdV) equations and Boussinesq models. Extensive studies on nonlinear dispersive equations have been conducted by Whitham (1974) and Zabusky (1965). These equations appear as asymptotic models. The presence of a quadratic nonlinearity in the wave speed introduces new challenges in stability and soliton interactions.

The most basic asymptotic dispersive equation is probably the nonlinear Schrodinger equation, which describes wave trains or frequency envelopes close to a given frequency, and their self-interactions (Taghizadeh and Mirzazadeh., 2012). Different methods have been used to obtain the solution of the Korteweg de Vries (KdV) equation. However, in this research work we construct a similar nonlinear equation of Schrodinger type with power nonlinearity($\sigma \geq 2$), wish to study its dynamic's nature and make comparisons with the existing nonlinear Schrodinger equation:

$$\begin{cases} i\partial_t u + \Delta u + \sigma |u|^p u = 0 \\ u(0,x) = u_0(x), x \in \mathbb{R}^n \end{cases}$$

In this study, we are interested in the dynamics' nature of the Cauchy problem with power-nonlinearity($\sigma \geq 2$).

$$\begin{cases} \psi_{tt} - [\alpha_1 + 3\alpha_2 \psi^2{}_x]\psi_{xx} + \alpha_3 \psi^\sigma = 0 & \text{in } \Omega \subset \mathbb{R}^+ \times \mathbb{R} \\ \psi(0,x) = \psi_0(x) \in S(\mathbb{R}) & \text{on } \partial\Omega \end{cases} \quad (1)$$

where $\psi: \Omega \subset \mathbb{R}^+ \times \mathbb{R} \to \mathbb{C}$ is a Schwartz function, $\alpha_j$ for $j = 1,2,3$ are arbitrary constants. This equation (1), without the power nonlinearity, is asymptotic model that showed up in studies like that of DNA dynamics for the descriptions of the dynamics of DNA molecules (Qian *et al*., 2014). The nonlinearity $\psi^\sigma$ here is used to capture the behaviour of nonlinear waves of optics, where wave collapse is possible, but interestingly in the setting where nonlinearity and dispersion behave according to the law $\psi_x^2 \psi_{xx}$.

This equation governs the evolution of nonlinear dispersive waves influenced by both quadratic nonlinearity and higher-order nonlinear terms. Such equations arise in nonlinear elasticity and nonlinear field theory. A similar structure is found in Boussinesq-type equations, which model shallow water waves, and in relativistic field equations with nonlinear potentials.

Furthermore, this equation, generalizes the standard Korteweg–de Vries (KdV) equation and Boussinesq equations, which have been extensively studied in soliton theory. Quadratic nonlinearity in the wave speed gives rise to a diverse range of wave phenomena, including soliton interactions, energy transfer mechanisms, and possible wave breaking scenarios.

Several works have investigated related nonlinear wave models. For instance, the Boussinesq equation for shallow water waves (Whitham, 1974) and the generalized KdV equation (Zabusky & Kruskal, 1965) have been analyzed using spectral and inverse scattering methods. More recent studies on nonlinear dispersive PDEs like Frauendiener et. al (2022), Sabo et. al. (2025). In the former, blow-up study is carried out on a 2- dimensional form of the NLS equation (1.0) in case $p = 2$, while in the latter, higher dispersive type of the equation is considered. For more available results and advanced studies on the dispersive equations read the monograph by Klien and Saut (2021). However, the study is motivated by the work of Bona & Smith (1975) and Craig et al. (2010), which explore stability and numerical approximations.

The nonlinear dispersive equations such as NLS equations, may have the tendencies of not admitting an exact solution, thus, in such cases numerical approach becomes the final resort. Since these equations are evolutionary, one finds a mapping method that uses a coordinates transformation that could transform them to a system of ordinary differential equations which are then time-integrated numerically. For instance, Niamh (2019) constructed an exact solutions for the NLS equation. Mehri (2013) presented the numerical Solutions of Korteweg de-Vries and Korteweg de Vries-Burger's equations using computer program.

**Methodology**

Let us consider some properties of the underlying equation. The equation (1), has, associated to it, the Hamiltonian density $\mathcal{H}$

$$\mathcal{H} = \tfrac{1}{2}[\psi^2_t + \left(\alpha_1 + \tfrac{\alpha_2}{2}\psi^2_x\right)\psi^2_x + 2\tfrac{\alpha_3}{\sigma+1}\psi^{\sigma+1} \tag{2}$$

the mass of the system is defined as:

$$M[\psi(t)] = \int_{\mathbb{R}} \psi(x,t)dx \tag{3}$$

The energy of the system is defined, using equation (2),

$$E[\psi(t)] = \int_{\mathbb{R}} \mathcal{H}dx$$
$$= \tfrac{1}{2}\int_{\mathbb{R}} \left(\psi^2_t - \left(\alpha_1 + \tfrac{\alpha_2}{2}\psi^2_x\right)\psi^2_x + \tfrac{2\alpha_3}{\sigma+1}\psi^{\sigma+1}\right)dx \tag{4}$$

Thus, mass (3) and energy (4) of the system are conserved if their time derivatives are zero, i.e.

$$\tfrac{d}{dt}M[\psi(t)] = 0 \text{ and } \tfrac{d}{dt}E[\psi(t)] = 0.$$

Both mass and energy are conserved quantities for the nonlinear wave equation under consideration, meaning that they remain constant in time during the evolution of the solution.

The underlying equation (1) is derived by using the variational derivative of the Hamiltonian $H(t) = E(t)$ and its density $\mathcal{H}$:

$$\frac{\delta H}{\delta \psi_t} = \frac{\partial \mathcal{H}}{\partial \psi_t} - \frac{\partial}{\partial t}\left[\frac{\partial \mathcal{H}}{\partial(\partial_t \psi_t)}\right] + \frac{\partial^2}{\partial t^2}\left[\frac{\partial \mathcal{H}}{\partial(\partial_{tt} \psi_t)}\right] + \cdots$$

whose Lagrangian density $\mathcal{L}$ is given by

$$\mathcal{L} = \Pi\psi_t - \mathcal{H} = \frac{1}{2}\psi^2_t + \frac{1}{2}\left(\alpha_1 + \frac{\alpha_2}{2}\psi^2_x\right)\psi^2_x - \frac{\alpha_3}{\sigma+1}\psi^{\sigma+1}$$

where $\Pi = \psi_t$.

**Dispersion relation**

Dispersion relation tells more how the wave propagates in space and time.

For small-amplitude waves, we assume a solution of the form

$$\psi(t,x) = e^{i(kx-\omega t)} \qquad (5)$$

Substituting (5) into the main equation (1), we obtain the dispersion relation:

$$-\omega^2 e^{i(kx-\omega t)} = \alpha_1(-ik)^2 e^{i(kx-\omega t)} + 3\alpha_2(-ik)^2 e^{i(kx-\omega t)} \cdot (-ik)^2 e^{2i(kx-\omega t)}.$$

Or

$$-\omega^2 \approx \alpha_1(-ik)^2 + 3\alpha_2(-ik)^4 \cdot [1 + O(kx-\omega)]$$

leading to the dispersion relation

$$\omega^2 \approx \alpha_1 k^2 - 3\alpha_2 k^4.$$

The dispersion relation is purely real if $\alpha_1 k^2 \geq 3\alpha_2 k^4$ and satisfies the relation that

$$v_g := \nabla\omega = \pm\frac{\alpha_1 k - 6\alpha_2 k^3}{\sqrt{\alpha_1 k^2 - 3\alpha_2 k^4}} \neq \frac{\omega}{k} := v_p.$$

meaning that the group velocity $v_g$ is not equal to the phase-velocity $v_p$. A classic example of dispersion relation is $\omega = \alpha k - \beta k^3$ for the Korteweg-de-Vries (KdV) equation:

$$\psi_t + \beta\psi\psi_x + \alpha\psi_{xxx} = 0.$$

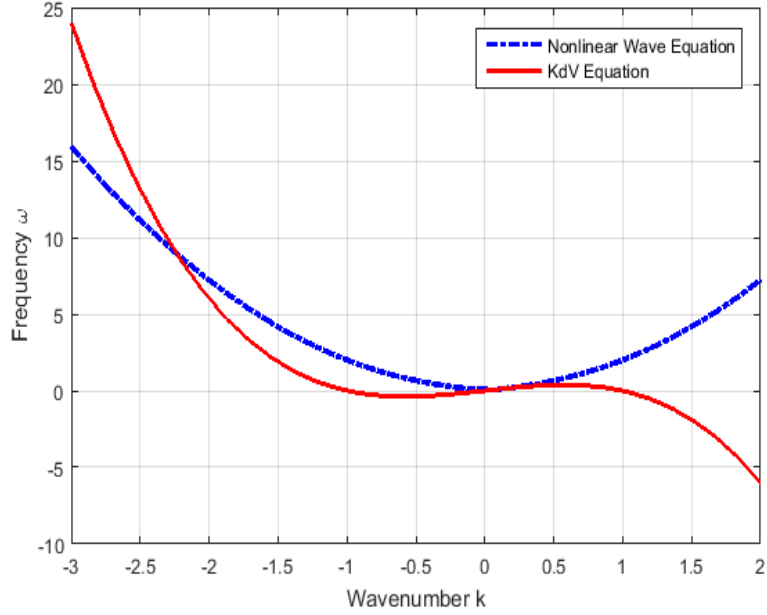

**Fig.1**: Dispersion relation for KdV and the nonlinear wave equations.

**Analytical & Spectral Approaches**

It is observed that the full linear form of the equation is the classical wave equation $\phi_{tt} = c^2 \phi_{xx}$, therefore, it is not an interesting problem worth studying at the moment. We therefore, consider dispersive-nonlinearity term while neglecting the power nonlinearity term $\psi^\sigma$.

Let us consider the partially-linearized form of the equation in (1):

$$\psi_{tt} - (\alpha_1 + \alpha_2 \psi_x^2)\psi_{xx} = 0.$$

We are interested in an exponentially decaying solution in the Schwartz Space $S(\mathbb{R})$, that is such $\psi$ having itself and its derivatives vanishing at the boundary. Using

$$\xi = x - ct, \quad \text{and} \quad \psi(x,t) \equiv \phi(\xi),$$

the equation transformed as

$$(c_1 - \alpha_1)\phi_{\xi\xi} - \alpha_2 \phi_\xi^2 \phi_{\xi\xi} = 0.$$

Upon integrating both-sides w.r.t. $\xi$ once, we have

$$(c^2 - \alpha_1)\phi_\xi - \frac{\alpha_2}{3}\phi_\xi^3 = C.$$

As the equation involves only first order derivative terms, we solve for $\phi_\xi$. However, given that

$$\lim_{\xi \to \pm\infty} \phi = \lim_{\xi \to \pm\infty} \phi' = \lim_{\xi \to \pm\infty} \phi'' = \cdots = 0,$$

then $C = 0$, the equation reduces to

$$(c^2 - \alpha_1)\phi_\xi - \frac{\alpha_2}{3}\phi_\xi^3 = 0.$$

The solution to the equation is

$$\phi_\xi = 0, \text{ or } \phi_\xi = \pm\sqrt{\frac{3}{\alpha_2}(c^2 - \alpha_1)}.$$

Implying that, for arbitrary constants, the solutions read:

$$\phi = c_1, \text{ or } \phi = \pm\sqrt{\frac{3}{\alpha_2}(c^2 - \alpha_1)}\xi + c_2.$$

It is observed that, the solution to the equation with dispersive-nonlinearity is

$$\psi = \text{const}, \quad \text{or} \quad \psi = \pm\sqrt{\frac{3}{\alpha_2}(c^2 - \alpha_1)}(x - ct) + \text{const}.$$

Obviously, this is not fully in the Schwartz space since $\psi$ itself is not in $S(\mathbb{R})$ but in the extended Schwartz Space $S_\omega(\mathbb{R})$, the *Beurling space*, that allows functions to grow at infinity with controlled by weight $\omega = \omega(x)$ while having the derivatives decaying sufficiently fast at $\infty$.

**Numerical Method: Spectral**

Spectral method is an efficient method that uses global basis and useful in determining the numerical solution of the of evolutionary PDEs. The method is applied as follows. As there are linear terms involved in the equation, we apply Fourier transform in space and use any finite difference techniques, such as RK4, to determine the time evolution of the equation.

**Fourier Transform**: Apply the Fourier transform of the both sides of the equation of the form

$$\psi_t = L\psi + N[\psi]$$

to have a set of ODEs of the form

$$\hat{\psi}_t = \mathcal{F}[L\psi] + \mathcal{F}[N[\psi]]$$

where $\mathcal{F}[u]$ stands for the Fourier transform of a function $u$ :

$$\mathcal{F}[u](k) = \int_{-\infty}^{\infty} e^{-ikx} u(x)dx := \hat{u}(k), \qquad \mathcal{F}^{-1}[\hat{u}(k)] = u(x) = \int_{-\infty}^{\infty} e^{ikx} \hat{u}(k)dk$$

and its inverse defined by $\mathcal{F}^{-1}[\hat{u}]$.

**Time Integration**: one use numerical time integration method, in our case, the RK4 method or simple time-stepping (e.g. Euler) as the case may be, to integrate the obtained ODE, in the Fourier Space.

**Inverse Fourier Transform**: invert the obtained solution for the ODEs in the Fourier space to the physical space. This way, one construct the solution to the equation.

**Implementation**: The given problem, is written as

$$\psi_t = \alpha_1 \psi_{xx} + [3\alpha_2 \psi_{xx} \psi_x^2 + \alpha_3 \psi^\sigma] = L\psi + N[\psi]$$

where the nonlinear term is $N[\psi] = 3\alpha_2 \psi_{xx} \psi_x^2 + \alpha_3 \psi^\sigma$ while the linear term is $\alpha_1 \psi_{xx}$.

The nonlinear part may be further split into two parts, the dispersive-nolinear and the power nonlinearity terms:

$$N[\psi] = D_N[\psi] + P_N[\psi].$$

In this regard, the Fourier transform of the equation gives

$$\hat{\psi}_t = \alpha_1(-ik)^2 \hat{\psi} + 3\alpha_2 \mathcal{F}[\psi_{xx} \psi_x^2] + \alpha_3 \mathcal{F}[\psi^\sigma].$$

Then the time stepping method is applied to solve the problem, numerically.

Let us describe how the Fourier transform is applied in solving the given PDE using the spectral method with RK4 time-stepping:

The given the fully nonlinear equation (1):

I. Convert to First-Order System of odes: introduce
$\hat{\psi} = \mathcal{F}[\psi], \quad V = \mathcal{F}[\psi_t],$
Then, rewrite as:
$$\hat{\psi}_t = V$$
$$V_t = \mathcal{F}\left[\alpha_1 + 3\alpha_2 \left(\mathcal{F}^{-1}[\hat{\psi}_x]\right)^2 \mathcal{F}^{-1}(\widehat{\psi_{xx}}) - \alpha_3 \left(\mathcal{F}^{-1}[\hat{\psi}]\right)^\sigma\right].$$
where $\hat{\psi}$ is the Fourier transform of $\psi$, and $V$ is the Fourier transform of $\psi$.

II. Spectral Differentiation: compute the $\psi_x$ and $\psi_{xx}$ in Fourier Space to have
$$\mathcal{F}[\psi_x] = (-ik)\,\hat{\psi}, \qquad \mathcal{F}[\psi_{xx}] = (-ik)^2\,\hat{\psi},$$
Then, re-compute the nonlinear terms in physical space:
$$\mathcal{F}^{-1}[\hat{\psi}] \mapsto \psi, \qquad \mathcal{F}^{-1}[\widehat{\psi_x}] \mapsto \psi_x, \qquad \mathcal{F}^{-1}[\widehat{\psi_{xx}}] \mapsto \psi_{xx};$$
Next, we compute
$$\mathcal{F}\left[\alpha_1 + 3\alpha_2 \left(\mathcal{F}^{-1}[\hat{\psi}_x]\right)^2 \mathcal{F}^{-1}(\widehat{\psi_{xx}})\right]$$
using Pseudo-spectral multiplication.

III. Time volution with RK4: Use the Runge-Kutta 4$^{th}$ order (RK4) method for time-stepping:
We define
- $k_1, k_2, k_3, k_4$ that approximate $\hat{\psi}$ at different stages.
- $l_1, l_2, l_3, l_4$ approximate $\hat{\psi}_t$ (denoted as $V$) at different RK4 steps.
These steps are computed as follows:
1. First step (Euler's method part):
$$k_1 = dt \cdot V;$$
$$l_1 = dt \cdot F(\widehat{\psi_t}, \quad V).$$
where $F(\widehat{\psi_t}, V)$ is the right-hand side of the equation.
2. Second step (Midpoint estimate 1):
$$k_2 = dt \cdot \left(V + \frac{1}{2}l_1\right);$$
$$l_2 = dt \cdot F\left(\hat{\psi} + \frac{1}{2}k_1, \quad V + \frac{1}{2}l_1\right).$$
3. Third step (Midpoint estimate 2):
$$k_3 = dt \cdot \left(V + \frac{1}{2}l_2\right);$$
$$l_3 = dt \cdot F\left(\hat{\psi} + \frac{1}{2}k_2, \quad V + \frac{1}{2}l_2\right)$$
4. Fourth step (Final correction):

$$k_4 = dt \cdot \left(V + \frac{1}{2}l_3\right);$$

$$l_4 = dt \cdot F\left(\hat{\psi} + \frac{1}{2}k_3, \quad V + \frac{1}{2}l_3\right).$$

5. Final update of $\hat{\psi}$ and $V$:

$$\widehat{\psi_t} = \hat{\psi}_t + \frac{1}{6}(k_1 + 2k_2 + 2k_3 + k_4);$$

$$V = V + \frac{1}{6}(\ell_1 + 2\ell_2 + 2\ell_3 + \ell_4),$$

where $V$ represents $\widehat{\psi_t}$, and $\ell_i$ are the intermediate approximations for its evolution.

Thus, $\ell_i$ for $i = 1,2,3,4$ are RK4 estimates for the acceleration $\widehat{\psi_{tt}}$ at different stages of the method.

**Results and Discussion**

The energy and mass of the system were shown to be conserved, confirming its Hamiltonian structure. Spectral methods combined with RK4 time integration were implemented to solve the equation efficiently. The stability of the numerical method was examined through mass and energy evolution, as well as the growth of Fourier modes. The interaction of solitons showed *elastic collisions* in some regimes, while energy was transferred in highly nonlinear regimes.

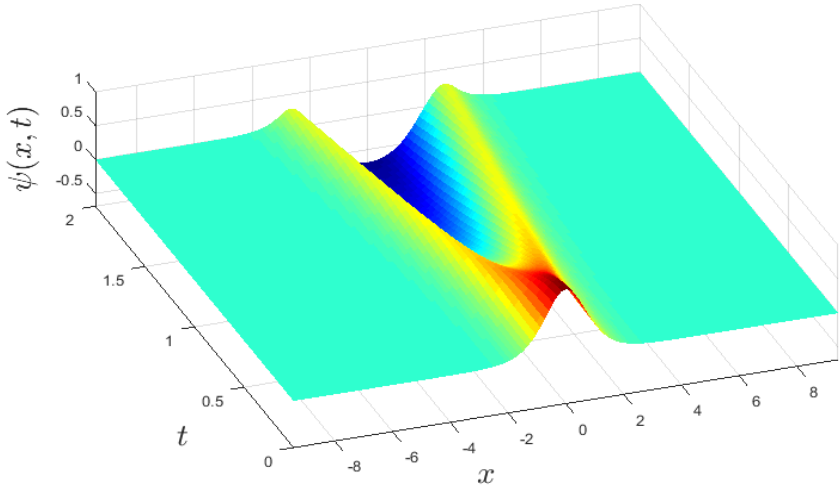

**Fig.1.** The close-up view of the solution $\psi$ for $\alpha_1 = \alpha_2 = \alpha_3 = 1$ with $\psi_0(x) = e^{-x^2}$.

These figures, Fig.1 and Fig.2., show the evolution of the solution $\psi(t, x)$. The solutions exhibit similar behaviour where a single initial soliton is decomposes into two solitons of virtually same size. They propagate over time and maintain their shapes. In either of the cases, defocusing (repulsive) nonlinearity ($\alpha_1 = \alpha_2 = -\alpha_3 = 1$) and attracting nonlinearity ($\alpha_1 = \alpha_2 = \alpha_3 = 1$) there is soliton fragmentation or fission behaviour as shown in Fig.1-2.

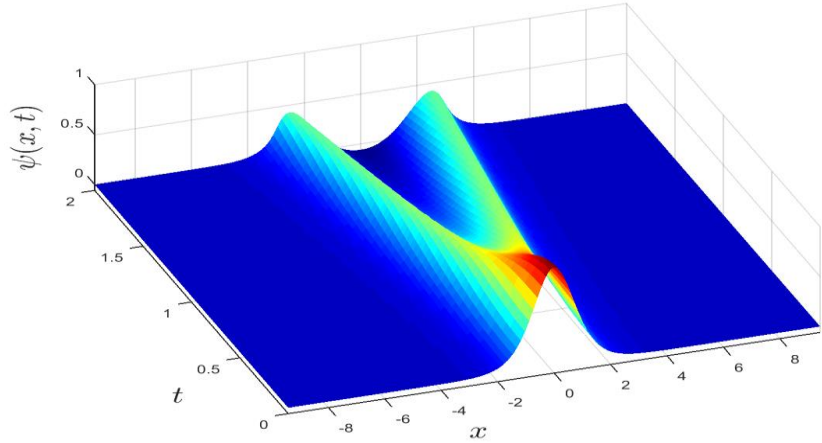

**Fig.2.** The close-up view of the solution $\psi$ for $\alpha_1 = \alpha_2 = -\alpha_3 = 1$ with $\psi_0(x) = e^{-x^2}$.

This phenomenon, of soliton fission, occurs as a result of energy redistribution that leads to the formation of multiple solitons from the initial solitary wave in nonlinear dispersive systems. The splitting is the effect of the dispersive nonlinearity term $\phi_x^2 \phi_{xx}$ .

Given that the equation is Hamiltonian, we could track the stability of the solution via evolution of the mass and energy. Moreover, to verify the numerical stability of the spectral method, one uses the Fourier spectrum which indicates rapid convergence when there is quick convergence of the Fourier modes.

In the result, we observed the evolution of the mass and the energy via the $\log_{10}-$plot as this helps to closely study the behaviour of the mass and energy at smaller scale. For the case when $\alpha_1 = \alpha_2 = \alpha_3 = 1$, the mass and energy evolution are presented in Fig.3.

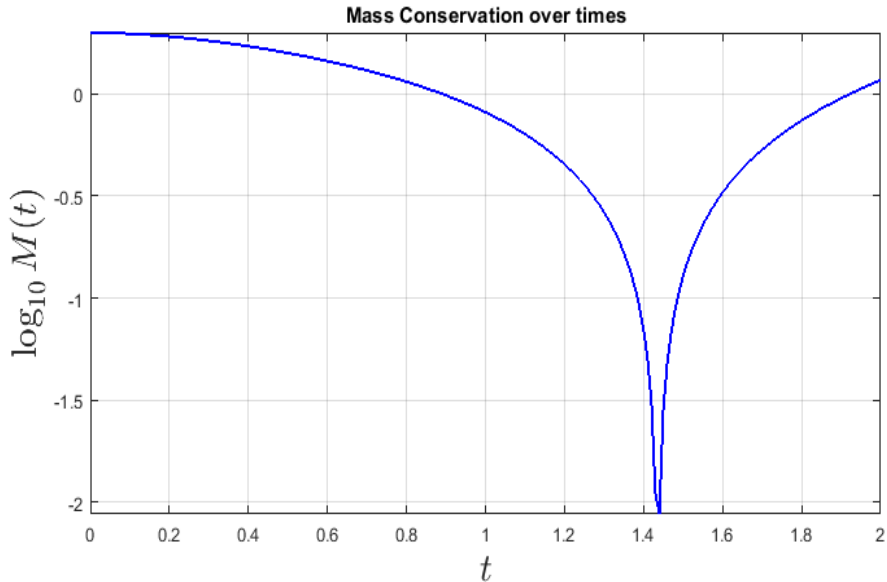

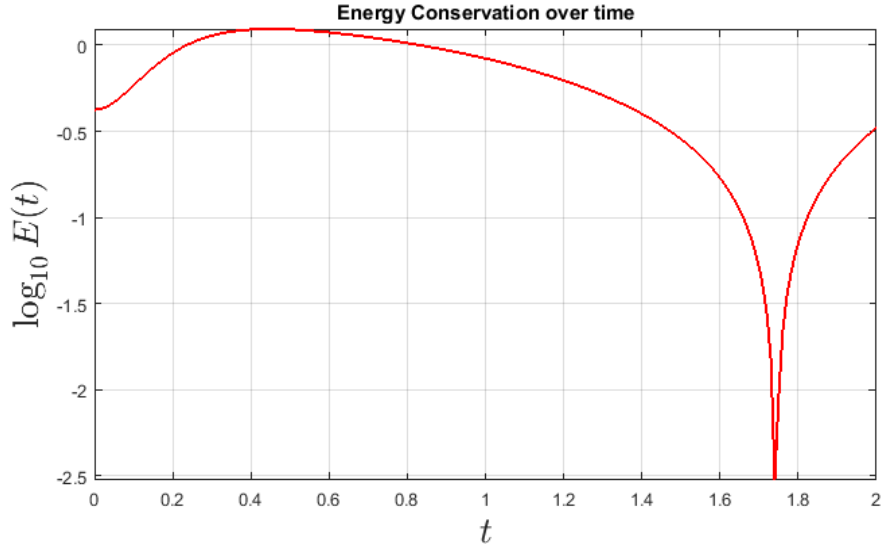

**Fig.3**. Evolution of the mass and energy functions.

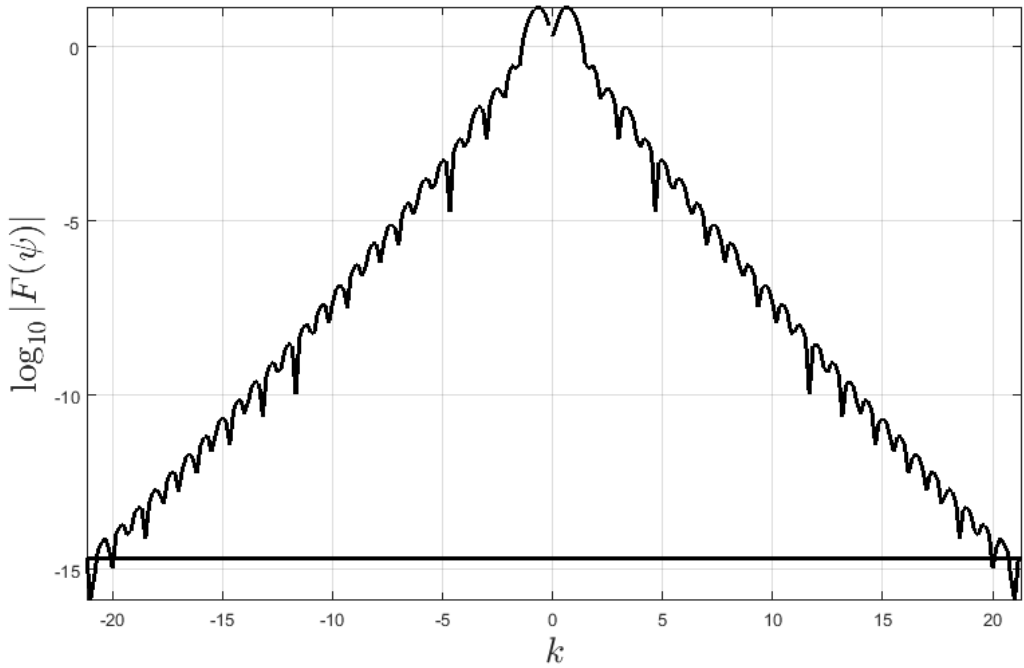

**Fig.4.** Fourier modes for the case $\alpha_1 = \alpha_2 = \alpha_3 = 1$, at the final time $t = 2$.

The Fourier spectrum in Fig. 4 shows remarkable fall-off up to $10^{-15}$, thereby showing stability of the spectral method. The evolutions of the mass and energy fluctuate and the effective behaviour are clearly understood for long time. However, longer time causes numerical instability of our method resulting from accumulation of numerical rounding errors.

On the other hand, it is important to check the stability of the solution to the equation per se. Therefore, we carry out the stability analysis of the solution associated with the equation (1) as follows. Since longer time could suffer potential numerical instability, it is wiser to study the stability in comparison to a perturbed solution. We will also keep track of the energy and the mass as they have to be bounded for the stability of the system.

We introduce small parameter $\varepsilon \ll 1$ and real number $k$ and assume the perturbation of the form

$$\psi_{\text{pert}}(t,x) = \psi(t,x) + \varepsilon \cos(k_p x) \qquad (6)$$

The growth-rate is computed via the expression

$$\gamma := \log_{10} \left| \frac{\|\psi_{\text{pert}}(t,x) - \psi(t,x)\|}{\|\psi(t,x)\|} \right| < \infty, \qquad (7)$$

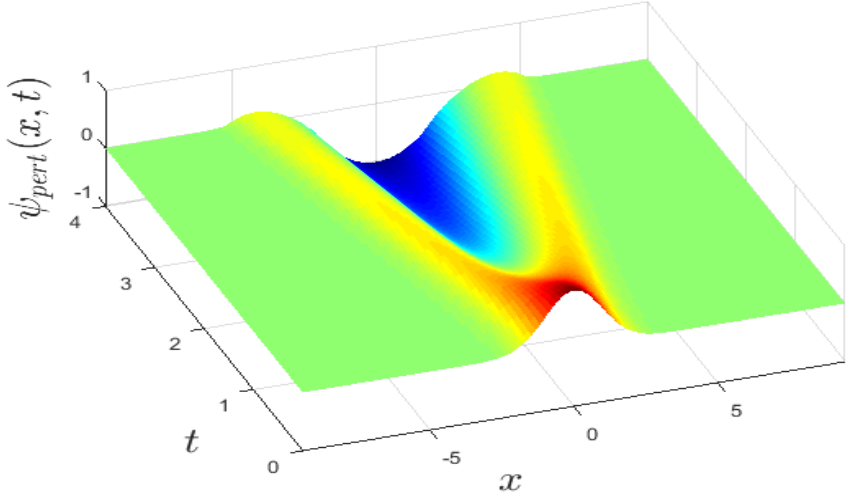

**Fig.5.** Plot of the perturbed solution with $\alpha_1 = \alpha_2 = \alpha_3 = 1$ and initial data $\text{sech}^2(x)$, with $\varepsilon = 10^{-3}$, and $k = 2$.

Fig.5. represents the plot of the perturbed solution which appears smaller in amplitude. In particular, we use $\varepsilon = 0.001$ with $\text{sech}(x)^2$ initial data and $k_p = 1/4$ in the perturbed term $\cos(k_p x)$. As this oscillates between $-1$ to $1$, then $\psi + \varepsilon \cdot \cos(k_p x)$, in (6), is not far away from the solution $\psi$ itself. The growth rate $\kappa$ in equation (7) is expected to be small.

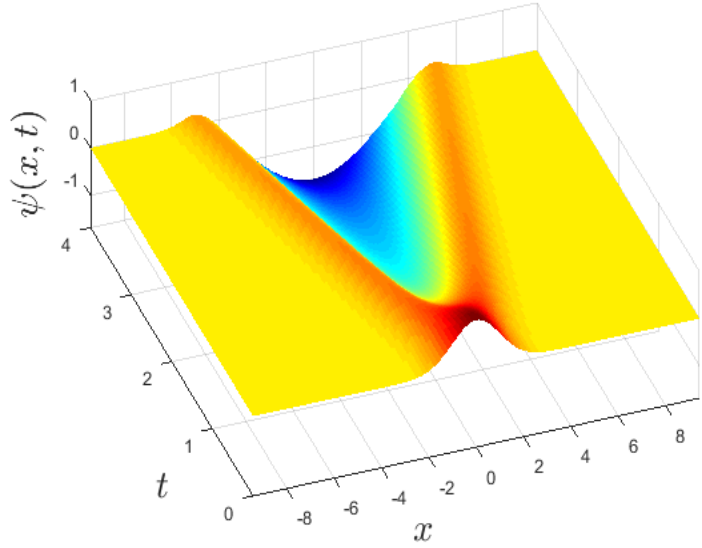

**Fig 6.** Plot of the unperturbed solution with $\alpha_1 = \alpha_2 = \alpha_3 = 1$ and initial data $\text{sech}^2(x)$.

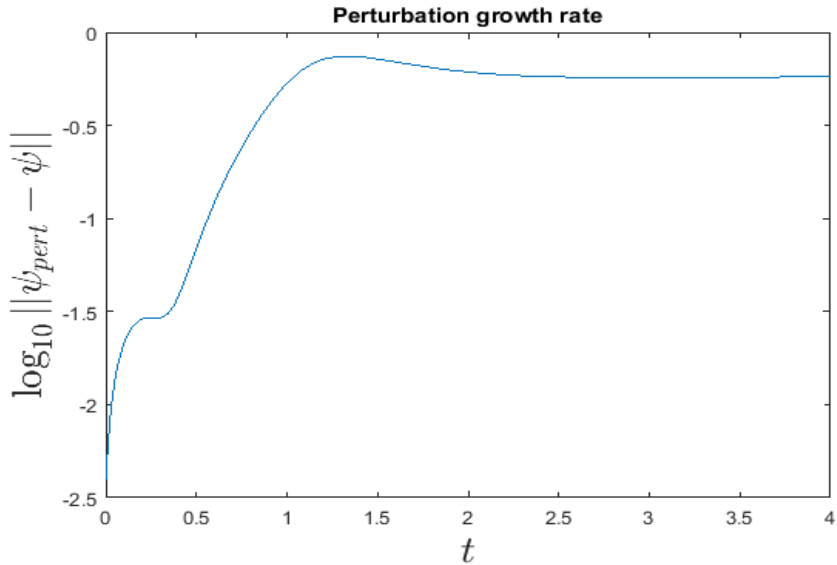

**Fig. 6.**: The growth rate via $\log$−plot of the difference $\psi_{\text{pert}} - \psi$ for longer time for $\alpha_1 = \alpha_2 = \alpha_3 = 1$ using initial data $\psi_0 = \text{sech}^2 x$.

The growth rate shows boundedness for longer time. In principle, for the stability of the system, we must have steady evolution of the $\gamma$ as $t \to \infty$. As it stays bounded, we conclude the underlying nonlinear equation (1) is stable. Similarly, we plot the stability of the spectral implementation in Fig.7. via the Fourier spectrum.

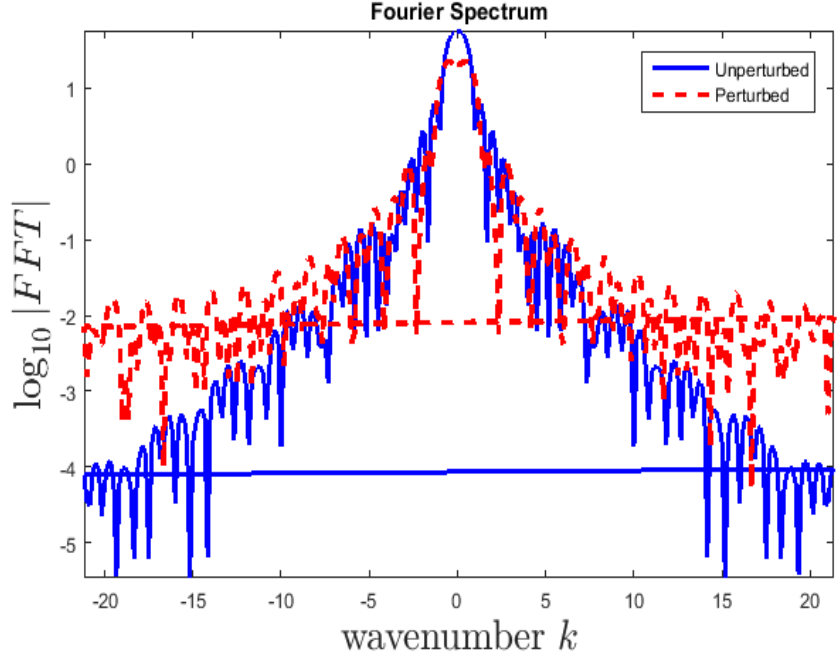

**Fig. 7.** The Fourier modes of the perturbed (solid line) and unperturbed (broken line) in the case $\alpha_1 = \alpha_2 = \alpha_3 = 1$ having initial data $\text{sech}^2(x)$.

Hence, the numerical results presented here suggested the following conjecture.

## 5. Main Result: Conjecture

I) For Schwartzian initial data $\psi_0(x)$ with positive coefficients $\alpha_1, \alpha_2, \alpha_2$ and $\sigma \geq 2$ the solution $\psi(t,x)$ to the equation (1) decomposes into small wave-packets.

II) For $\psi_{\text{pert}}(t,x)$, then

$$\gamma := \log_{10}\left[\frac{\|\psi_{\text{pert}}(t,x) - \psi(t,x)\|}{\|\psi(t,x)\|}\right] < \infty,$$

provided that $\alpha_1, \alpha_2, \alpha_2$ stay positive, where $\|\psi\|$ is the norm of $\psi$.

III) $\lim_{t \to \infty} \|\psi(t,x)\|$ stays finite for small $t$ if $\alpha_1, \alpha_2, \alpha_2$ are negative values for any $\sigma \geq 2$.

This conjecture implies that (I) a Schwartzian initial lump decomposes into smaller lumps, any other possible combination of different signs would lead to singularity of the solution (II) the stability is guaranteed over $t$ provided that the quantity $\gamma < \infty$ stays finite over time; and (III) Solution to the equation exists locally, i.e. for finite time. This is to say that singularity or blow-up sets in after a finite time. In fact, the singularity emerges at any time $t \geq \frac{1}{2}$.

**Conclusion**

This study investigated a nonlinear dispersive wave equation with quadratic nonlinearity in the wave speed. The spectral method was implemented to solve the PDE numerically, demonstrating energy and mass conservation and whose stability is determined via its Fourier spectrum. The dispersion relation was compared with the classical KdV equation, highlighting key differences in wave propagation behaviour. We show that, via conjecture, that the equation is stable in the focusing case and blow-up is possible for repulsive nonlinearity due to the hyperbolic nature of the equation. Future work could focus on establishing the proof of the stated conjecture on the stability of the solution, application of alternate numerical approaches and other integrability properties of the equation can be explored.